\documentclass[onefignum,onetabnum]{siamart220329}

\usepackage{lipsum}
\usepackage{amsfonts}
\usepackage{graphicx}
\usepackage{epstopdf}
\usepackage{algorithmic}
\ifpdf
  \DeclareGraphicsExtensions{.eps,.pdf,.png,.jpg}
\else
  \DeclareGraphicsExtensions{.eps}
\fi

\newsiamremark{remark}{Remark}
\newsiamremark{hypothesis}{Hypothesis}
\crefname{hypothesis}{Hypothesis}{Hypotheses}
\newsiamthm{claim}{Claim}

\headers{Weak Distance between 2-sphere and Approximations}{K. Koga}

\title{Computing weak distance between the 2-sphere and its nonsmooth approximations\thanks{
\funding{This work was partially supported by, JST ACT-X Grant Number JPMJAX2106, JST CREST Grant Number JPMJCR22P5, and JSPS KAKENHI Grant Number\,21K20325.}}}

\author{Kazuki Koga\thanks{\email{koga.k.ac@m.titech.ac.jp}.}}

\usepackage{amsopn}

\ifpdf
\hypersetup{
  pdftitle={Computing weak distance between the 2-sphere and its nonsmooth approximations},
  pdfauthor={K. Koga}
}
\fi

\DeclareFontFamily{U}{mathx}{}
\DeclareFontShape{U}{mathx}{m}{n}{<-> mathx10}{}
\DeclareSymbolFont{mathx}{U}{mathx}{m}{n}
\DeclareMathAccent{\widecheck}{0}{mathx}{"71}

\begin{document}

\maketitle

\begin{abstract}
A novel algorithm is proposed for quantitative comparisons between compact surfaces embedded in the three-dimensional Euclidian space. The key idea is to identify those objects with the associated surface measures and compute a weak distance between them using the Fourier transform on the ambient space. In particular, the inhomogeneous Sobolev norm of negative order for a difference between two surface measures is evaluated via the Plancherel theorem, which amounts to approximating a weighted integral norm of smooth data on the frequency space. This approach allows several advantages including high accuracy due to fast-converging numerical quadrature rules, acceleration by the nonuniform fast Fourier transform, and parallelization on many-core processors. In numerical experiments, the 2-sphere, which is an example whose Fourier transform is explicitly known, is compared with its icosahedral discretization, and it is observed that the piecewise linear approximations converge to the smooth object at the quadratic rate up to small truncation.
\end{abstract}




\section{Introduction}\label{sec1}
Quantitative comparisons between two-dimensional surfaces play crucial roles in geometry processing tasks such as shape classifications and reconstructions \cite{JiWaHuMaNi2019,IcDeTu2023}, and their importance is further increasing due to disciplines of machine learning that manipulates geometric data. In those applications, the so-called Chamfer distance, which measures a difference between two point clouds by averaging the Euclidean distance from elements of one point cloud to the other and vice versa, is a rough but popular way for comparing a surface to another. However, naive sampling on geometric objects ignores their ``curved" structures, and hence the Chamfer distance is often unsatisfactory in scientific computing involving discretization of surfaces. In fact, computationally feasible comparisons between the 2-sphere $\mathbb{S}^2$ and its piecewise linear approximations are still a nontirivial problem in spite of their potential applications to convergence study of mesh generation \cite{BaFr1985,WaLe2011} and justification of remeshing \cite{KoHu2020, McCoBu2018}, both of which require at least a few digits of accuracy in an appropriate distance that takes the two-dimensional structure into account. \par
From the viewpoint of mathematics, there are several concepts that quantify differences between geometric objects and define convergence in suitable senses. Among those, the Hausdorff distance \cite{BrHa2013} on metric spaces is an intuitive tool and extensively studied in Riemannian geometry, whereas its numerical implementation is not realistic except for simple cases such as small point clouds. On the other hand, recent developments in computational optimal transport enable to approximate the Wasserstein distance between probability measures on the Euclidian space \cite{Villani2009, PeCu2019}. In the context of comparisons between two-dimensional surfaces, relevant contributions include, but are not limited to, the cases for two point clouds with help of an entropic regularization \cite{Cuturi2013}, for a set of simplices and a point cloud based on the damped Newton's method \cite{MeMeTh2018}, and for two disk-like surfaces via conformal maps \cite{LiDa2011}. Nevertheless, it is also difficult to extend the techniques to the situation where both objects are two-dimensional, and the major barrier to computations of the above concepts is that they require global optimization involving the background metric. \par
In the planar setting, Koga \cite{Koga2021} introduces an invariant-based approach to comparing a smooth periodic curve to another via the Fourier coefficients of the arclength variable, which are accessible from arbitrary parametrizations and preserved under changes of parametrizations up to translation and reflection, and it allows to evaluate quality of nonuniform discretization of periodic planar curves with local complexity \cite{Koga2022}. However, this approach heavily relies on the simple topology of periodic curves and their smoothness. Therefore, we immediately encounter difficulties when extending the idea to smooth surfaces having the topology of the $2$-sphere $\mathbb{S}^2$, which do not admit an analogy of the Fourier coefficients of the arclength variable, and to comparisons between smooth and piecewise smooth objects, where the decay rates of their invariant Fourier transforms are drastically different.\par
To devise a new comparison algorithm that covers a wider range of applications, the present work switches the measure-theoretic viewpoint of optimal transport to that of harmonic analysis. Namely, compact surfaces embedded in the three-dimensional Euclidian space $\mathbb{R}^3$ are identified with the associated surface measures, and their Fourier transforms are computed with the structures of the ambient space. Since the Paley-Wiener theorem guarantees smoothness of the Fourier transforms of compactly supported surface measures, this strategy translates the problem into comparisons between well-behaved functions on the frequency space and allows to compute the inhomogeneous Sobolev norm of negative order via the Plancherel theorem. Such a real-analytic concept should be weaker than the Hausdorff and even the Wasserstein distances because its definition does not directly use the metric on $\mathbb{R}^3$.
\par From computational perspectives, however, the framework of Fourier analysis leads to several advantages including accuracy, efficiency, and parallelizability. First, the Fourier transforms of surface measures are reduced to integrals on canonical domains via changes of variables and evaluated by fast-converging numerical quadrature rules there, while weighted integral norms of the resulting smooth functions are well approximated by the three-dimensional trapezoidal rule over bounded cubes. Second, the application of quadrature rules on canonical domains corresponds to generating weighted Dirac measures on $\mathbb{R}^3$ as good approximations of surface measures, and the resulting point clouds and data on the trapezoidal grid in the frequency space are efficiently bridged by the nonuniform fast Fourier transformm (NUFFT) \cite{BaMaKl2019}. Last, the linearity of the Fourier transform enables to decompose a compact surface into a set of disjoint pieces and apply quadrature rules in an element-wise manner. This feature is a basis for parallelization on massively parallel architectures called Graphics Processing Units (GPUs). In numerical experiments, the 2-sphere in $\mathbb{R}^3$, an example whose surface measure has an explicit form of its Fourier transform in terms of the Bessel function, is compared with the icosahedral discretization \cite{BaFr1985}, and it is observed that the simplex-based approximations converge in the Sobolev norm to the smooth object at the quadratic rate up to small truncation. Furthermore, we claim that the suggested algorithm outperforms an existing method for computing the Fourier transform of 2-simplices based on the exact formula \cite{IcDeTu2023}.\par
The rest of this paper is organized as follows. Section \ref{formulate} formulates the comparison problem for compact surfaces via the Fourier transform on $\mathbb{R}^3$ and introduces relevant examples of surface measures. Section \ref{method} explains numerical methods and some techniques in their efficient implementations. Section \ref{results} presents numerical results and evaluates the implemented codes. Concluding remarks are given in Section \ref{conclusion}.


\section{Formulations}\label{formulate}
This section provides basic concepts in Fourier analysis on $\mathbb{R}^3$ with a few examples of Borel measures whose Fourier transforms are explicitly known. For conciseness, the content is focused on how to formulate quantitative comparison between compact surfaces in the framework of Fourier analysis, and readers interested in more details are referred to standard textbooks on Euclidian harmonic analysis \cite{Grafakos2014-1,Grafakos2014-2}.
%
%
%
%
%
%
%
\subsection{Fourier transform on $\mathbb{R}^3$}
Given a Borel-measurable function $f$ on $\mathbb{R}^3$, its Fourier transform $\widehat{f}$ is defined as the integrals
\begin{equation}
\label{eq:ft_forward}
\widehat{f}(\xi) = \int_{\mathbb{R}^3} f(x) e^{-2\pi i \xi \cdot x} dx,\quad \xi \in \mathbb{R}^3,
\end{equation}
where $dx$ denotes the Lebesgue measure and we call the variable $\xi$ the frequency. In the setting of the whole space $\mathbb{R}^3$, one of the most important examples is the heat kernel
\begin{equation}
\label{eq:hear_kernel}
\phi(t,x) = \frac{1}{(4\pi t)^{3/2}} e^{-|x|^2/4t},\quad x \in \mathbb{R}^3,\quad t>0,
\end{equation}
which has the Fourier transform of the form
\begin{equation}
\label{eq:heat_ft}
\widehat{\phi}(t,\xi) = e^{-4\pi^2 t |\xi|^2},\quad \xi \in \mathbb{R}^3,\quad  t>0.
\end{equation} 
For a wide class of functions, the inversion formula
\begin{equation}
\label{eq:ft_inverse}
\widecheck{f}(x) = \int_{\mathbb{R}^3} f(\xi) e^{2\pi i \xi \cdot x} d\xi,\quad x \in \mathbb{R}^3,
\end{equation}
serves as the inverse operation of (\ref{eq:ft_forward}) in a suitable sense. For instance, the heat kernel (\ref{eq:hear_kernel}) is obtained in the pointwise sense by applying the operation (\ref{eq:ft_inverse}) to its Fourier transform  (\ref{eq:heat_ft}).\par
If $f$ belongs to the class $L^2(\mathbb{R}^3)$ of square-integrable functions, the norm $\|f\|_{L^2(\mathbb{R}^3)}$ is directly computed from its Fourier transform $\widehat{f}$ via the Plancherel theorem
\begin{equation}
\label{eq:plancherel}
\|f\|^2_{L^2(\mathbb{R}^3)}=\int_{\mathbb{R}^3} |f(x)|^2 dx=\int_{\mathbb{R}^3} |\widehat{f}(\xi)|^2 d\xi,
\end{equation}
which says that the Fourier transform is an isometry on $L^2(\mathbb{R}^3)$. Moreover, the Fourier transform can be used to define several function spaces that takes regularity of functions into account. Among those,  the so-called (inhomogeneous) Sobolev space, denoted by $W_s^p(\mathbb{R}^3)\, (1<p<\infty,\,s\in \mathbb{R})$, is the set of Borel-measurable functions for which the integral norm
\begin{equation}
\label{eq:sobolev_p}
\| f\|_{W_s^p(\mathbb{R}^3)} =\|((1+|\,\cdot\,|^2)^\frac{s}{2} \widehat{f})\,\widecheck{ }\,\|_{L^p(\mathbb{R}^3)},
\end{equation}  
is finite. In the following, we call the exponent $s$ the order of regularity, and for $s>0$ multiplying the weight $(1+|\xi|^2)^\frac{s}{2}$ amounts to taking derivatives of the function $f$, while for $s<0$ it acts as smoothing on $f$. Then, by the relation (\ref{eq:plancherel}), the norm of the Sobolev space for the case $p=2$, denoted by $H_{s}(\mathbb{R}^3)$, is equivalently defined as the integral on the frequency space
\begin{equation}
\label{eq:sobolev_2}
\| f\|_{H_{s}(\mathbb{R}^3)} =\biggl( \int_{\mathbb{R}^3} (1+|\xi|^2)^s |\widehat{f}(\xi)|^2 d\xi\biggr)^\frac{1}{2},
\end{equation}
and it enables to compute the norm $\| f\|_{H_s(\mathbb{R}^3)}$ without the inversion formula (\ref{eq:ft_inverse}).
 \par
As an extension of (\ref{eq:ft_forward}), the Fourier transform of a Borel measure $\mu$ is defined as
\begin{equation}
\label{eq:ft_measure}
\widehat{\mu}(\xi) = \int_{\mathbb{R}^3}  e^{-2\pi i \xi \cdot x} d\mu(x),\quad \xi \in \mathbb{R}^3,
\end{equation}
but the inversion formula (\ref{eq:ft_inverse}) is no longer used to recover $\mu$. For instance, the Dirac measure $\delta_0$ at the origin is Borel and its Fourier transform is
\begin{equation}
\label{eq:ft_dirac}
\widehat{\delta_0}(\xi) = e^{-2\pi i \xi \cdot 0}=1,\quad \xi \in \mathbb{R}^3,
\end{equation}
which is smooth but only belongs to the class of bounded functions. However, the Sobolev norm (\ref{eq:sobolev_2}) can be meaningful for such measures with $s<0$, which is equivalent to computing the $L^2$ norm after smoothing via the convolution with a well-known function called the Bessel kernel $G_s$ defined by $\widehat{G_s}(\xi)=(1+|\xi|^2)^\frac{s}{2}$ \cite{Mattila2015}. Furthermore, such a procedure is extended to elements of the space of tempered distributions $\mathcal{S}'(\mathbb{R}^3)$, which is the space of continuous linear functionals on the Schwartz space $\mathcal{S}(\mathbb{R}^3)$. In fact, the set of all finite Borel measures can be regarded as a subset of $\mathcal{S}'(\mathbb{R}^3)$, and thus it is possible to define distances between them as functionals in terms of the Sobolev norm with an appropriate value of $s$.
%
%
%
%
%
%
%
\subsection{Surface measures}
Among Borel measures on $\mathbb{R}^3$, the special class called surface measures plays a crucial role in the present work. Roughly speaking, the surface measure of a two-dimensional surface in $\mathbb{R}^3$ defines integrals of measurable functions with respect to its area element. For example, the surface measure of the $2$-sphere $\mathbb{S}^2 = \{x\in\mathbb{R}^3\,|\, |x|=1\}$, denoted by $\sigma$, gives rise to integrals
\begin{equation}
\label{eq:def_surf_measure_sphere}
\int_{\mathbb{R}^3} f(x) d\sigma(x) = \int_{\mathbb{S}^{2}}  f(x) dS,
\end{equation}
where $dS$ is the area element of $\mathbb{S}^2$, and its Fourier transform is explicitly written:
\begin{equation}
\label{eq:ft_nsphere}
\widehat{\sigma}(\xi) = \int_{\mathbb{S}^{2}}  e^{-2\pi i \xi \cdot x} dS = \frac{2\pi}{|\xi|^{\frac{1}{2}}}J_{\frac{1}{2}}(2\pi|\xi|).
\end{equation}
Here, the function $J_\nu$ is the Bessel function of order $\nu$, and using the gamma function $\Gamma$, it is conveniently represented by the series expansion
\begin{equation}
\label{eq:Bessel_series}
J_\nu (t) = \sum^\infty_{j=0} \frac{(-1)^j}{j!} \frac{1}{\Gamma(j+\nu+1)}\biggl(\frac{t}{2} \biggr)^{2j+\nu},
\end{equation}
which implies that the Fourier transform (\ref{eq:ft_nsphere}) is a power series whose radius of convergence is $\infty$. Thus, the Fourier transform $\widehat{\sigma}$ is smooth on $\mathbb{R}^3$, and due to the asymptotic expansion \cite{Stein1993}
\begin{equation}
\label{eq:asymp_bessel}
J_\nu(t) \sim \biggl(\frac{2}{\pi}\biggr)^{\frac{1}{2}} t^{-\frac{1}{2}}\cos\biggl(t- \frac{\pi\nu}{2} - \frac{\pi}{4}\biggr),
\end{equation}
it has a slow decay rate $\widehat{\sigma}(\xi)=\mathcal{O} (|\xi|^{-1})\,(|\xi|\rightarrow \infty)$. These properties are also found by simplifying the formula (\ref{eq:ft_nsphere}) with the identities (\ref{eq:Bessel_series}) and $\Gamma(\frac{1}{2}+n) = \frac{(2n-1)!!}{2^n}\sqrt{\pi}$:  
\begin{equation}
\label{eq:ft_nsphere_sinc}
\widehat{\sigma}(\xi) =\frac{2\sin(2\pi|\xi|)}{|\xi|}.
\end{equation} 
In general, a wide class of measures $\mu$ supported on hypersurfaces in $\mathbb{R}^d$ with nonzero Gaussian curvature satisfy 
\begin{equation}
\label{eq:ft_measure_hypersurf}
|\widehat{\mu}(\xi)| \leq C(\mu) (1+|\xi|)^{-\frac{d-1}{2}},
\end{equation}
where the constant $C$ depends on both the support of $\mu$ and the distribution on it \cite{Demeter2020}. Since the surface $\mathbb{S}^2$ is written as the union of the northern and southern hemispheres whose overlap is a zero set with respect to its surface measure, the estimate (\ref{eq:ft_measure_hypersurf}) also applies to (\ref{eq:ft_nsphere}) and is regarded as a generalization of its decay property. Moreover, finiteness of the Sobolev norm of Borel measures in $\mathbb{R}^d$ is examined via the so-called $\alpha$-energy \cite{Mattila2015}, which is defined as 
\begin{equation}
\label{eq:def_alpha_energy}
I_\alpha(\mu) = \iint \frac{1}{|x-y|^\alpha}d\mu(x)d\mu(y)= c(d,\alpha) \int_{\mathbb{R}^3}  |\xi|^{\alpha-d}|\widehat{\mu}(\xi)|^2d\xi,
\end{equation} 
where $0<\alpha<d$ and $c(d,\alpha)$ is a constant determined by $d$ and $\alpha$. Since the dimension of surface measures, which is equal to $2$ for surfaces in $\mathbb{R}^3$, can be also understood as the supremum of the parameters $\alpha$ for which the integral (\ref{eq:def_alpha_energy}) is finite, it is clear that surface measures of interest in the present work belong to the Sobolev space $H_s(\mathbb{R}^3)$ for $s<-\frac{1}{2}$. On the other hand, it is well-known that the Fourier transforms of finite Borel measures of compact support are smooth and extended to entire functions. This fact is often referred to as the Paley-Wiener theorem \cite{Hormander2015}, and allows us to assume that compact surfaces in $\mathbb{R}^3$ have the associated surface measures whose Fourier transforms are smooth and design numerical quadrature rules based on the favorable property.    \par
Another important example of surface measures is those of affine simplices, which are fundamental for approximating smooth objects. Let $Q_0^k$ be the standard $k$-simplex
\begin{equation}
\label{eq:simplex_st}
Q^k_0 = \{(x_1,\ldots,x_k)\in \mathbb{R}^k\,|\,x_1+\cdots +x_k \leq 1,\, x_i\geq 0 \,(i=0,\ldots,k)\}.
\end{equation}
For example, the simplex $Q_0^2$ is the triangle with the vertices at $0, e_1, e_2$, and $Q_0^3$ is the tetrahedron with $0,\,e_1,\,e_2,\,e_3$, where $0$ is the origin and $\{e_i\}^k_{i=1}$ is the standard basis of $\mathbb{R}^k$. Then, an affine $k$-simplex in $\mathbb{R}^d\,(d\geq k)$ with vertices $\{p_0,p_1,\ldots, p_k\}$, denoted by $Q^k$, is the image of the standard simplex $Q_0^k$ by the affine map
\begin{equation}
\label{def_affine_simplex}
L(x) = p_0 + Ax,\quad x\in Q^k_0,
\end{equation}
where the $d\times k$ matrix $A$ satisfies $Ae_i = p_i - p_0\,(i=1,\ldots,k)$ \cite{Rudin1976}. Note that an affine simplex is normally defined as an oriented object, but the order of the set $\{p_0,p_1,\ldots, p_k\}$ is not considered here because our interest lies in integrals with respect to its surface measure. The Fourier transforms of the surface measures of affine simplices also have closed forms for general $d$ and $k$ \cite{JiWaHuMaNi2019}, and the present work is in particular related to the case for $d=3$ and $k=2$. By an abuse of notation, the surface measure of a $k$-simplex $Q^k$ is denoted by the same symbol, and its Fourier transform is written as
\begin{equation}
\label{eq:simplex32_ft}
\widehat{Q^k}(\xi) = -\frac{\gamma}{4\pi^2} \sum^2_{i=0}\frac{e^{-2\pi i\xi\cdot p_i}}{\prod_{j\neq i}\{\xi\cdot (p_i-p_j)\} },\quad \xi \in \mathbb{R}^3,
\end{equation}
where $\gamma$ is the ratio of the surface area of $Q^k$ to that of $Q^k_0$. Again, the function (\ref{eq:simplex32_ft}) is smooth due to the compactness of $Q^k$ and the Paley-Wiener theorem, whereas its decay property is somewhat different from that of the sphere $\mathbb{S}^2$. In fact, writing $\xi = a\xi_N +\xi_\perp$, where $\xi_N$ is a unit normal vector of $Q^k$ and $\xi_\perp$ is orthogonal to $\xi_N$, it is easy to see that, for example,
\begin{equation}
\label{eq:simplex32_ft_orth}
\widehat{Q^k}(\xi) = e^{-2\pi i a\xi_N\cdot p_0}\widehat{Q^k}(\xi_\perp) ,
\end{equation}
and the absolute value $|\widehat{Q^k}(\xi) |$ is constant in the parameter $a$. As a surface measure in $\mathbb{R}^3$, however, it also belongs to the space $H_s(\mathbb{R}^3)$ for $s<-\frac{1}{2}$, and thus comparisons between the $2$-sphere $\mathbb{S}^2$ and its approximations by affine simplices can be formulated as measuring the difference between (\ref{eq:ft_nsphere}) and a linear combination of (\ref{eq:simplex32_ft}) in the Sobolev norm (\ref{eq:sobolev_2}).

\section{Numerical methods}\label{method}
This section is devoted to describing numerical methods for computing the distance between the 2-sphere $\mathbb{S}^2$ and its simplex-based approximation in the sense of the Sobolev norm (\ref{eq:sobolev_2}). The main algorithm consists of numerical quadrature rules on simplices and the frequency space, respectively, and a fast summation algorithm that bridges those schemes in an efficient way.
%
%
%
%
%
%
%
\subsection{Quadrature rules on simplex}
For computing the Sobolev norm (\ref{eq:sobolev_2}), the first step is to approximate the Fourier transforms of given measures to sufficiently large frequencies. Due to the linearity of the Fourier transform, this task for a simplicial complex is reduced to computing the Fourier transform of each simplex and summing up their contributions in the frequency space. In fact, several previous studies including \cite{JiWaHuMaNi2019,IcDeTu2023} employ the exact formula (\ref{eq:simplex32_ft}) and implement the direct summation in the frequency space. However, the complexity of such native computations is $\mathcal{O}(N_S\cdot N_F)$, where $N_S$ is the number of simplices and $N_F$ the number of target points on the frequency space. Besides, the closed form (\ref{eq:simplex32_ft}) suffers from large cancellation errors when a simplex is almost degenerate or when the vector $\xi$ is almost orthogonal to one of the edges of a simplex, and practical implementations circumvent this difficulty by adding small regularization terms to the denominator \cite{JiWaHuMaNi2019} or replacing the formula with the case for $p_1=p_2$ or $p_2=p_3$ \cite{IcDeTu2023}, both of which limit accuracy for  data representing small-scale structures.  \par
In the present work, we choose a more general approach to approximating the Fourier transforms of $2$-simplices using a combination of a simple change of variables and fast-converging numerical quadrature rules on the standard simplex (\ref{eq:simplex_st}). Namely, with the affine map $L$ given by (\ref{def_affine_simplex}), perform the change of variables
\begin{equation}
\label{eq:change_variable}
\int_{Q^k} f(x) dS = \int_{Q_0^k} f(L(u)) \gamma_A du,
\end{equation}
where the constant $\gamma_A$ is the Jacobian determined by the matrix $A$, and a quadrature rule $\{(u_i,q_i)\}$ on $Q^k_0$ is applied to the right-hand side of (\ref{eq:change_variable}): 
\begin{equation}
\label{eq:change_variable_approx}
 \int_{Q_0^k} f(L(u)) \gamma_A du \approx \sum_{i}  f(L(u_i)) q_i \cdot \gamma_A.
\end{equation}
Since the nodes $\{L(u_i)\}$ and the weights $\{q_i \cdot \gamma_A\}$ are independent of the integrand $f$,  this procedure is equivalent to converting the quadrature rule $\{(u_i,q_i)\}$ on $Q^k_0$ to  $\{(L(u_i),q_i \cdot \gamma_A)\}$ on $Q^k$. In other words, this approach can be regarded as a construction of approximate measures in terms of the Dirac measures on $\mathbb{R}^3$ from numerical quadrature rules on the parameter domain, and it turns out that approximating a surface measure by weighted atomic ones allows to use efficient algorithms for summing up their contributions on the frequency space. In our implementation, the 8-point Gauss-Legendre rule is applied twice to the right-hand side of (\ref{eq:change_variable}) as an iterated integral and it serves as the quadrature rule $\{(u_i,q_i)\}$ on the standard domain $Q^k_0$. We are aware that more symmetric rules are available on triangles \cite{ZhCuLi2009, WaXi2003}, which requires to find symmetric polynomials and perform numerical optimization while an iteration-free algorithm is available for the Gauss-Legendre rules \cite{Bogaert2014}.
%
%
%
%
%
%
%
\subsection{Quadrature rules on $\mathbb{R}^3$}
Once adequate methods are found for obtaining data in the frequency space, the next step is to design numerical quadrature rules on the whole $\mathbb{R}^3$ that approximate the Sobolev norm (\ref{eq:sobolev_2}). Since integrands are evaluated at finite points only, it is first necessary to choose a bounded domain $D$ where the truncated integral is sufficiently close to that on $\mathbb{R}^3$.  To this end, we employ the cube $D=[-\xi_\text{max},\xi_\text{max}]^3$ and sample integrands on the lattice
\begin{equation}
\label{eq:def_trap_grid}
\xi_{i,j,k} = (ih,jh,kh),\quad h=\frac{2\xi_\text{max}}{M}, \quad i,j,k=-\frac{M}{2},\ldots, \frac{M}{2}.
\end{equation}
As we see later, this specific choice is compatible with fast evaluations of the Fourier transform if the positions of atomic measures as approximations of surface measures are appropriately scaled. Then, the three-dimensional trapezoidal rule is applied to the integral over $D$:
\begin{equation}
\label{eq:trap_rule_R3}
\int_D f(\xi) d\xi\approx  \sum^{\frac{M}{2}}_{i=-\frac{M}{2}}\sum^{\frac{M}{2}}_{j=-\frac{M}{2}}\sum^{\frac{M}{2}}_{k=-\frac{M}{2}} f(\xi_{i,j,k}) w_{i,j,k} h^3,
\end{equation}
where $w_{i,j,k}$ is 1 for points in the interior of the cube, $1/2$ for the faces, $1/4$ for the edges, and $1/8$ for the corners. For sufficiently large $\xi_\text{max}$, which depends on the decay property of a given integrand, the formula (\ref{eq:trap_rule_R3}) also approximates the integral over the entire $\mathbb{R}^3$. On discretization errors, we refer to the Euler–Maclaurin formula \cite{TrWe2014}, which explains the asymptotic errors of the one-dimensional trapezoidal rule for a smooth function $f$ in terms of its derivatives of odd orders $f^{(2n-1)}$:
\begin{equation}
\label{eq:euler-maclaurin}
I_M - I \sim \sum^\infty_{n=1} h^{2n} \frac{B_{2n}}{2n!}[f^{(2n-1)}(b) - f^{(2n-1)}(a)].
\end{equation}
Here, $I$ is the integral of $f$ over a closed interval $[a,b]$, $I_M$ is its approximation by the $M$-point trapezoidal rule, and $B_k$ is the $k$-th Bernoulli number. Thus, the order of the trapezoidal rule is $\mathcal{O}(h^2)$ for general $f$, but actual errors are expected to be rather small if the derivatives $f^{(2n-1)}$ decay rapidly. This observation is a basis for applications of the trapezoidal rule to fast-decaying functions including the heat kernel (\ref{eq:hear_kernel}), and a similar situation is found for the Sobolev norm (\ref{eq:sobolev_2}) if the Fourier transform of a given measure is smooth, as for finite Borel measures of compact support, and the absolute value of the negative order $s$ is large. However, we remark that the trapezoidal rule is not necessarily the best method for approximating this kind of integrals. For instance, due to the radial symmetry of the weight in (\ref{eq:sobolev_2}), one may consider integrals over balls centered at the origin and construct numerical quadrature rules on them.
%
%
%
%
%
%
%
\subsection{Nonuniform fast Fourier transforms}
As above, numerical quadrature rules are chosen for both affine simplices and $\mathbb{R}^3$ as the frequency space, but we still lack an efficient way of computing data on the latter from those on the former. By construction, the direct summation of contributions from $N_P=64\,\cdot\,N_S$ source points to $N_F=(M+1)^3$ targets on the frequency space has $\mathcal{O}(N_P\cdot N_F)$ complexity, which is essentially the same as the method based on the exact formula (\ref{eq:simplex32_ft}). If source points are also evenly spaced on a lattice such as (\ref{eq:def_trap_grid}), it is possible to perform this summation in $\mathcal{O}(N_P\log N_P)$ complexity using the so-called fast Fourier transform (FFT).  However, the approximation (\ref{eq:change_variable_approx}) leads to highly nonuniform point distributions on the affine simplex, and therefore is not compatible with the standard FFT. Instead, an efficient algorithm called the (Type-1) nonuniform fast Fourier transform (NUFFT) is capable of approximating the following sum to a prescribed relative accuracy $\epsilon_\text{rel}$: 
\begin{alignat}{2}
\label{eq:def_nufft_type1}\widehat{f}_\eta=\sum_{j=1}^{N_P}f_je^{-i\eta \cdot x_j}, \quad & \eta=(\eta_1, \eta_2, \eta_3)\in \mathbb{Z}^3 \cap [-\eta_\text{max}, \eta_\text{max}]^3,
\end{alignat}
where the coefficients $\{f_j\}$ are any complex numbers and the coordinates $\{x_j\}$ lie in $[-\pi,\pi)^3$. To utilize the standard FFT, the Type-1 NUFFT regards the exponential sum (\ref{eq:def_nufft_type1}) as the exact Fourier transform of weighted Dirac measures $ \mu = \sum_{j=1}^{N_P}f_j\delta_{x_j}$, and replaces it with
\begin{equation*}
\label{def_nufft_conv}
 \mu * \psi = \sum_{j=1}^{N_P}f_j \psi(\cdot-x_j), \quad \psi: \mbox{Gaussian, Kaiser-Bessel, etc.},
\end{equation*}
where $*$ is the convolution on the torus $\mathbb{T}^3=[-\pi,\pi)^3$ and $\psi$ is a smooth periodic function with a compact numerical support (e.g., Gaussian and Kaiser-Bessel kernels). Then, the Fourier transform of the regularized measure is computed by the three-dimensional trapezoidal rule with the FFT acceleration, and dividing the results by the symbols $\widehat{\psi}(\eta)$ yields approximations of the desired sums. The Type-2 NUFFT, which interpolates the Fourier series at nonuniform points, is also defined as the adjoint operation of the Type-1 NUFFT. As opposed to the direct summation, the complexity of this approximation algorithm is $\mathcal{O}(N_F\log N_F + N_P \log^3 (1/\epsilon_\text{rel}))$, and the number of source points $N_P$ and that of target points $N_F$ are no longer related to each other. Moreover, the convergence of such algorithms are analyzed for several kernels \cite{DuRo1993, BaMaKl2019, KeKuPo2009}, and now there are general-purpose implementations with their unique advantages \cite{BaMaKl2019, KeKuPo2009, GrLe2004}. Among those, FINUFFT \cite{BaMaKl2019}, which chooses $\psi$ to be the ``exponential of semicircle" kernel with nearly optimal aliasing errors \cite{Barnett2021}, is currently the fastest library with the OpenMP parallelization and is used in our implementation. Although the trapezoidal rule on $\mathbb{R}^3$ requires integrand values on $h\mathbb{Z}^3$ rather than $\mathbb{Z}^3$, by rewriting the summation of interest as
\begin{alignat}{2}
\label{eq:source_scale}
\sum_{j=1}^{N_P}f_je^{-ih\eta \cdot x_j}= \sum_{j=1}^{N_P}f_je^{-i\eta \cdot  hx_j},
\end{alignat}
we find that the evaluations on $h\mathbb{Z}^3$ are reduced to the Type-1 NUFFT with the coordinates $\{hx_j\}$. That is, for sufficiently large $M$, the Fourier transform on the trapezoidal grid (\ref{eq:def_trap_grid}) is obtained from a more clustered point cloud in the vicinity of the origin. Needless to say, this technique is peculiar to the truncated trapezoidal rule on $\mathbb{R}^3$, and in general it is required to compute contributions from nonuniform sources to nonuniform targets, which can be accelerated by the Type-3 NUFFT \cite{LeGr2005} with larger computational costs than by the Type-1 NUFFT . \par
As a relevant previous study, we conclude this section by mentioning the algorithm called the geometric NUFFT by Sammis and Strain \cite{SaSt2009} and its variant by Strain \cite{Strain2018}. The strategy in \cite{SaSt2009} replaces the Dirac measures in the Type-1 NUFFT above with polynomial distributions supported on simplices and directly convolutes the special class of measures with the B-spline kernel by precomputing quadrature rules on simplices that are exact up to some degree of the Bernstein polynomials. Since those exact quadrature rules also approximate smooth distributions in general, the geometric NUFFT can be regarded as a general-purpose algorithm for measures supported on simplicial complexes. However, it is difficult to extend such an approach to general smooth surfaces because an adequate polynomial basis cannot be found there and constructions of quadrature rules with exactness are impossible, whereas the parametrization-based approach in the present work is applicable to more general cases where surfaces are represented by smooth maps from bounded domains. 


\section{Numerical results}\label{results}
Now we show numerical results to validate the algorithm described in Section \ref{method}. In the following, the implemented CPU codes are run on MacBook Pro 2021 with Apple M1 Max (10 cores) and 32 GB memory, and the relative accuracy for the Type-1 NUFFT algorithm is set to $\epsilon_\text{rel}=10^{-15}$ (i.e. double precision) unless otherwise stated.
%
%
%
%
%
%
%
\begin{figure}[t]
\begin{center}
 \includegraphics[width=\linewidth]{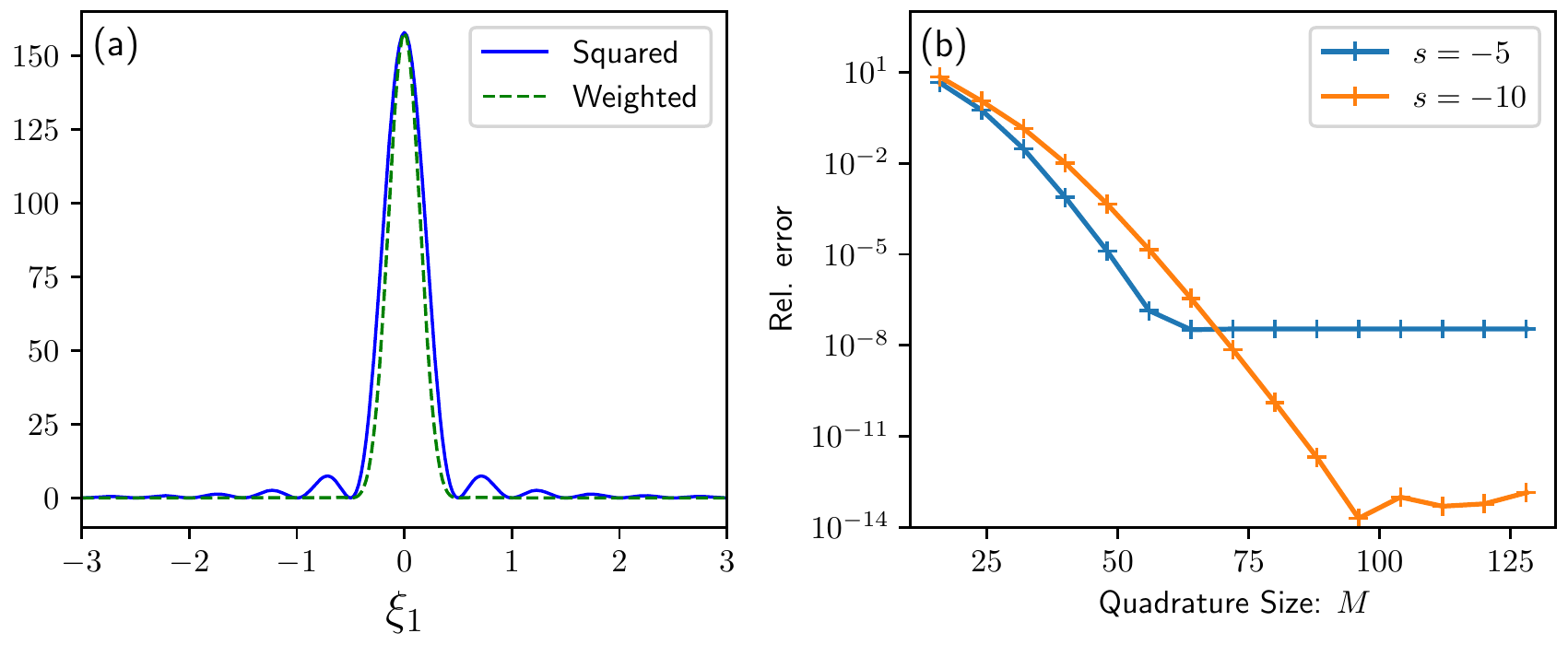}
\end{center}
\caption{\label{fig:plot_bessel_trap} Convergence of trapezoidal rule (\ref{eq:trap_rule_R3}) on $D=[-\xi_\text{max},\xi_\text{max}]^3$ : (a) exact  $|\widehat{\sigma}|^2$ and weighted function (\ref{eq:ft_2sphere_quad1D}) for $s=-10$, and (b) plots of relative errors for $s=-5,-10$ and $\xi_\text{max}=5$. } 
\end{figure}
\subsection{Convergence of trapezoidal rule}
We start numerical experiments by checking accuracy of the trapezoidal rule on the frequency space. As an example, we employ the exact representation (\ref{eq:ft_nsphere_sinc}) and consider the function
\begin{equation}
\label{eq:ft_2sphere_quad1D}
 g_s(\xi) =(1+|\xi|^2)^s|\widehat{\sigma}(\xi)|^2,\quad \xi \in \mathbb{R}^3,
\end{equation}
which is the integrand of the Sobolev norm $\|\sigma\|_{H_s(\mathbb{R}^3)}$ having the closed form
\begin{equation}
\label{eq:sphere_sobolev_exact}
\|\sigma\|_{H_s(\mathbb{R}^3)}^2 = \frac{4\pi^{\frac{3}{2}}}{\Gamma(-s)} \biggl(\Gamma\biggl(-s-\frac{1}{2}\biggr) - 2(2\pi)^{-s-\frac{1}{2}}K_{-s-\frac{1}{2}}(4\pi)\biggr),
\end{equation}
where $s<-\frac{1}{2}$ and $K_\nu$ is the modified Bessel function of the second kind of order $\nu$. This result is calculated by simplifying the integral of (\ref{eq:ft_2sphere_quad1D}) under the spherical coordinates and exploiting a convenient representation of the Bessel kernel $G_s$ in terms of the special function $K_\nu$ \cite{ArSm1961}.
Due to the radial symmetry, Figure \ref{fig:plot_bessel_trap}(a) shows the functions $|\widehat{\sigma}|^2$ (solid line) and $g_{s}$ for $s=-10$ (dashed line) on the $\xi_1$-axis with the notation $\xi=(\xi_1, \xi_2, \xi_3)$.  Since the function $|\widehat{\sigma}|^2$ decays as just $\mathcal{O}(|\xi|^{-2})$, it still exhibits visible oscillations away from the origin $\xi_1=0$ and hence the trapezoidal rule on the cube $D=[-\xi_\text{max},\xi_\text{max}]^3$ is practically of second order. On the other hand, the weighted function $g_{-10}$ decays as $\mathcal{O}(|\xi|^{-22})$ and the oscillatory behavior of the sinc function is no longer visible, while most energy around the origin is yet maintained. \par 
To see actual convergence of the trapezoidal rule for the function $g_s$, Figure \ref{fig:plot_bessel_trap}(b) plots relative errors $|I-I_M|/I$ for $s=-5,-10$ and $\xi_\text{max}=5$ versus the quadrature size $M$, where $I$ and $I_M$ are the exact value (\ref{eq:sphere_sobolev_exact}) and its approximations by the formula (\ref{eq:trap_rule_R3}), respectively. 
 Although both error curves show rather slow decays for small $M$, which is due to the large diameter of the integral domain and the oscillations of the integrand, the convergence of the trapezoidal rule is clearly faster than the theoretical order $\mathcal{O} (h^2)$. However, the crucial difference between the two cases is that even for intermediate $M$ relative errors for $s=-5$ are dominated by the truncation error at the level of $10^{-7}$, which also arises from the size of the domain $D$, while the total error continues to decrease for $s=-10$ and eventually reaches the level of $10^{-13}$. Thus, the trapezoidal rule (\ref{eq:trap_rule_R3}) efficiently approximates the integral of the weighted function (\ref{eq:ft_2sphere_quad1D}) for $s=-10$ over $D=[-5,5]^3$ and simultaneously gives a good estimate on the integral (\ref{eq:sphere_sobolev_exact}). Nevertheless, we again remark that efficiency of the trapezoidal rule largely depends on the pair $(s,\xi_\text{max})$, and to achieve a given rate of convergence and/or a prescribed truncation error, the maximal frequency $\xi_\text{max}$ must be prohibitively large as the order $s$ approaches the critical value $s^*=-\frac{1}{2}$.
%
%
%
%
%
%
%
\begin{figure}[t]
\begin{center}
 \includegraphics[width=\linewidth]{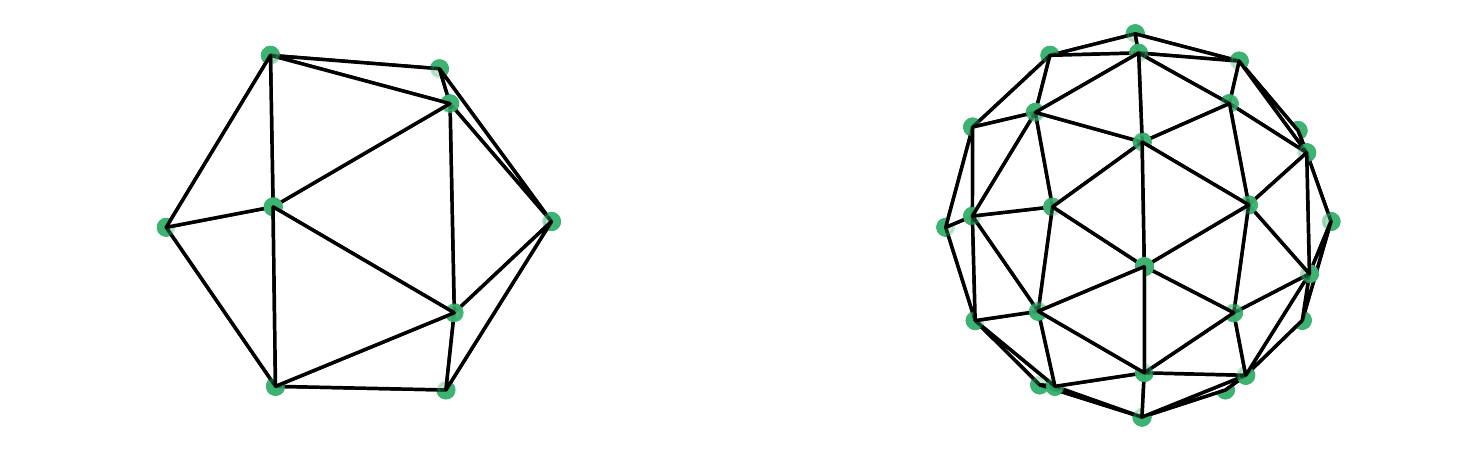}
\end{center}
\caption{\label{fig:plot_isocad} Illustrations of icosahedral discretization: Initial icosahedron $\mathcal{T}_0$ (left), and polyhedron $\mathcal{T}_1$ after first subdivision (right).} 
\end{figure}\par
 \subsection{Icosahedral discretization}
Next, we assemble the separate components in Section \ref{method} into a single algorithm that performs comparisons between the $2$-sphere $\mathbb{S}^2$ and its approximations. For this purpose, we introduce a recursive process called the icosahedral discretization \cite{BaFr1985}, which is one of the most intuitive methods for constructing simplex-based approximations of $\mathbb{S}^2$ and yet practical in scientific computing \cite{BoKeKrJa2017}. Figure\,\ref{fig:plot_isocad} shows two polyhedra generated from the icosahedral discretization. First, the process begins with an icosahedron inscribed in $\mathbb{S}^2$, denoted by $\mathcal{T}_0$, which is a regular polyhedron with the largest number of vertices uniformly distributed on $\mathbb{S}^2$. Since it is no longer possible to increase the number of vertices in a uniform manner, we consider a refinement of each face of the initial icosahedron and project new vertices onto $\mathbb{S}^2$. To this end, an additional vertice is chosen to be the midpoint of each edge and connected to the others on the same triangle. This step divides each face into four sub-triangles and the projection of those new vertices onto $\mathbb{S}^2$ leads to a finer approximation $\mathcal{T}_1$ of the smooth surface. Then, if the $(N-1)$-th refinement $\mathcal{T}_{N-1}$ of the initial icosahedron is obtained, one can repeat the subdivision of triangles again and generate the $N$-th refinement $\mathcal{T}_{N}$ in the same way as for the case $N=1$. This recursive process yields $10\times4^N + 2$ vertices and $20\times 4^N$ triangles for $\mathcal{T}_N$, and the order of the icosahedral discretization is expected to be $\mathcal{O}(4^{-N})$ in some distances as a two-dimensional analogue of the linear interpolation on closed intervals.\par

\begin{figure}[t]
\begin{center}
 \includegraphics[width=\linewidth]{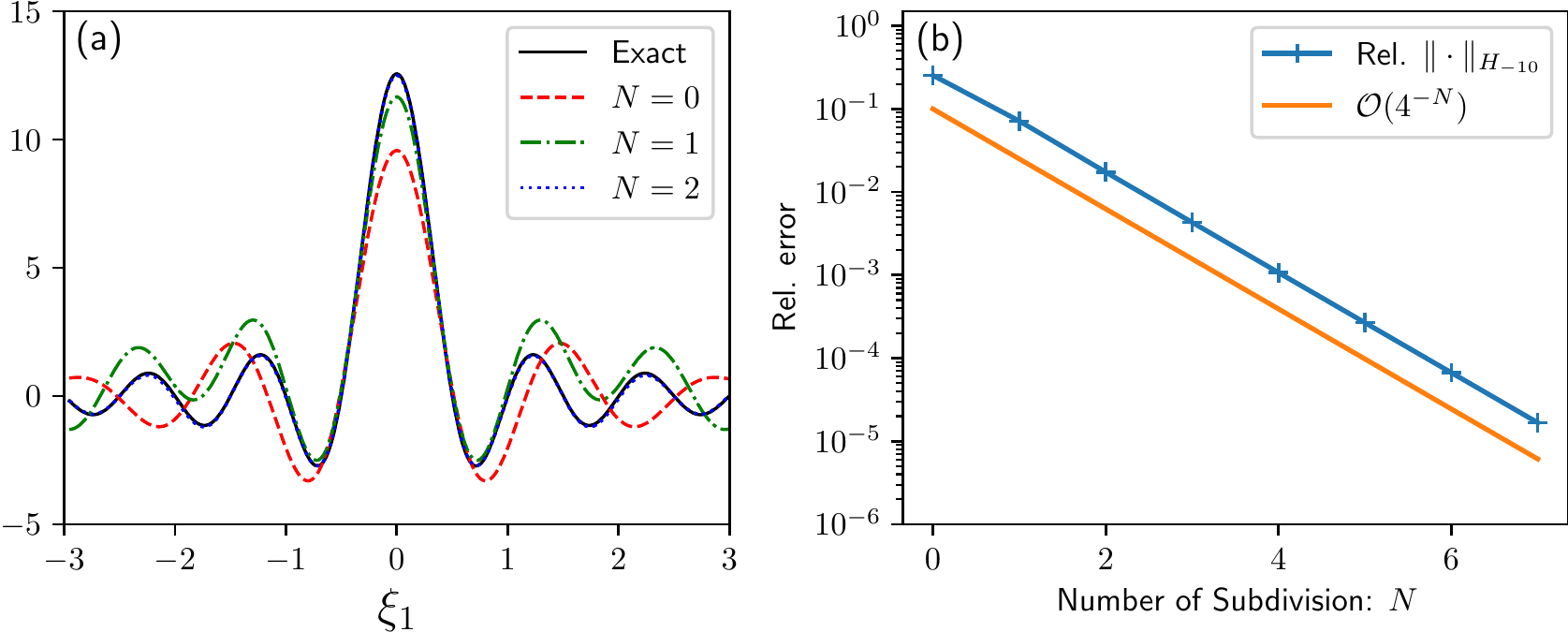}
\end{center}
\caption{\label{fig:conv_isocad}Convergence of icosahedral discretization : (a) plots for exact $\widehat{\sigma}$ and $\widehat{\mathcal{T}_n}$ with $N=0,1,2$, and (b) quadratic convergence in $\|\,\cdot\,\|_{H_{-10}(\mathbb{R}^3)}$ and theoretical $\mathcal{O}(4^{-N})$. }
\end{figure}

As numerical tests, we first observe behavior of simplex-based objects $\{\mathcal{T}_N\}$ from the icosahedral discretization in the frequency space, and then compute discretization errors in the sense of the Sobolev norm (\ref{eq:sobolev_2}) with the suggested algorithm. In the following, the icosahedral discretization is performed using the existing MATLAB code called Spherepts developed by Wright and Michaels\,\cite{WrMc2018}, which has a concise interface for generating several kinds of discretization of $\mathbb{S}^2$. Figure\,\ref{fig:conv_isocad}(a) illustrates pointwise comparisons on the $\xi_1$-axis between the Fourier transform of the surface measure $\sigma$, which is evaluated from the exact formula (\ref{eq:ft_nsphere_sinc}), and the real part of that of the polyhedron $\mathcal{T}_N$ for a few values of $N$. As easily seen, increasing $N$ leads to improved approximations of the function $\widehat{\sigma}$ (denoted by ``Exact") and schematically the graph for $N=2$ is indistinguishable from the true curve. Here, we remark that relative errors from the pointwise comparisons above are expected to become worse as the value of $\max\{|\xi_1|,|\xi_2|,|\xi_3|\}$ increases, and therefore the convergence is not necessarily uniform on $\mathbb{R}^3$. Next, Figure\,\ref{fig:conv_isocad}(b) shows relative errors in terms of the Sobolev norm (\ref{eq:sobolev_2}) for $s=-10$  versus the number of subdivisions $N$ with a theoretical $\mathcal{O}(4^{-N})$ curve. For the trapezoidal rule, we use $\xi_\text{max}=5$ and the grid size $M=128$, and it is clear that the error curve is almost parallel to $\mathcal{O}(4^{-N})$ and reaches the level of $10^{-5}$ for $N=7$.  Again, this rate of convergence is speculated as a two-dimensional analogue of the linear interpolation on closed intervals, and we conjecture that it derives from quadtratic pointwise convergence on the whole frequency space. However, quadratic convergence in the Sobolev norm is not necessarily deduced from such convergence in the pointwise sense, and we conclude that quadratic convergence in the Sobolev norm is numerically observed up to small truncation.
 \subsection{Performance evaluations} 
Here, we discuss further acceleration of the suggested algorithm with help of Graphics Processing Units (GPUs), which are massively parallel architecture with fast memory access and underly the success of recent machine learning applications, and evaluate performance of the implemented CPU/GPU codes for various $N$ and $M$. As seen in Section \ref{method}, the suggested algorithm converts quadrature rules on the standard simplex $Q^2_0$ to those on affine $2$-simplices in $\mathbb{R}^3$ and these operations can be performed in a simplex-wise manner. Moreover, a simplex-based approximation of a curved surface can consist of more than 100\,000 triangles and generating the corresponding point cloud multiplies the number by the size of a quadrature rule on $Q^2_0$. Such a large task as a set of tiny independent computations perfectly fits in GPU environments, because the devices achieve efficiency with their high memory bandwidth for supplying data to thousands of cores, and we parallelize those linear algebraic operations by the CUDA C/C++ language\,\cite{ChGrMc2014}.    \par 
As opposed to the previous step, the NUFFT algorihthm is more structured and allows less room for parallelization on GPUs. Since the fastest GPU version of the uniform FFT is provided by CUDA Toolkit\,\cite{ChGrMc2014}, most effort for accelerating the Type-1 NUFFT is devoted to the convolution part of the algorithm. In this context, there are currently two general-purpose GPU implementations called CUNFFT\,\cite{KuKu2012} and cuFINUFFT\,\cite{ShWrAnBlBa2021}, and the present work employs the latter that cleverly exploits a hardware resource unique to the CUDA environment. Roughly speaking, the core ideas of cuFINUFFT are (1) to sort a set of nonuniform points for improved memory access patterns, and (2) to decompose the whole convolution into subproblems on the so-called shared memory, which is essentially a programmable L1 cache and hence small but significantly faster than the global memory (i.e. GPU memory), for mitigating the issue of atomic collisions.  Since an initial indexing of nonuniform points may be highly random in the computational domain $\mathbb{T}^3=[-\pi,\pi)^3$, the sorting step reorders the input data so that two nonuniform points having consecutive indices contribute to close points on the regular grids or equivalently close memory addresses. On the other hand, accumulating contributions in parallel from numerous sources onto a single array causes the issue of atomic collisions, where memory access by  GPU threads is blocked until the ``current" thread completes its load and store. This inefficiency should become a bottleneck for parallel NUFFT algorithms if nonuniform points as input data are highly clustered in some small regions, and such situations are found in the use of the Type-1 NUFFT for evaluating the Fourier transform on the trapezoidal grid via the relation (\ref{eq:source_scale}). To reduce the chance of atomic collisions, the cuFINUFFT code decomposes the sorted point set into smaller subsets and adds their contributions to subgrids on the shared memory, and those intermediate results are separately stored on the global memory for obtaining the final results of the entire convolution later. One drawback of their technique is that the shared memory is too small for practical computations in the double-precision floating-point arithmetics. Nevertheless, it is found that the single-precision Type-1 NUFFT is accurate enough for reproducing the results in Figure\,\ref{fig:conv_isocad} and is expected to be practical in various applications. More details of cuFINUFFT including benchmarks against other implementations are thoroughly examined by its authors \cite{ShWrAnBlBa2021}. \par
\begin{figure}[t]
\begin{center}
 \includegraphics[width=\linewidth]{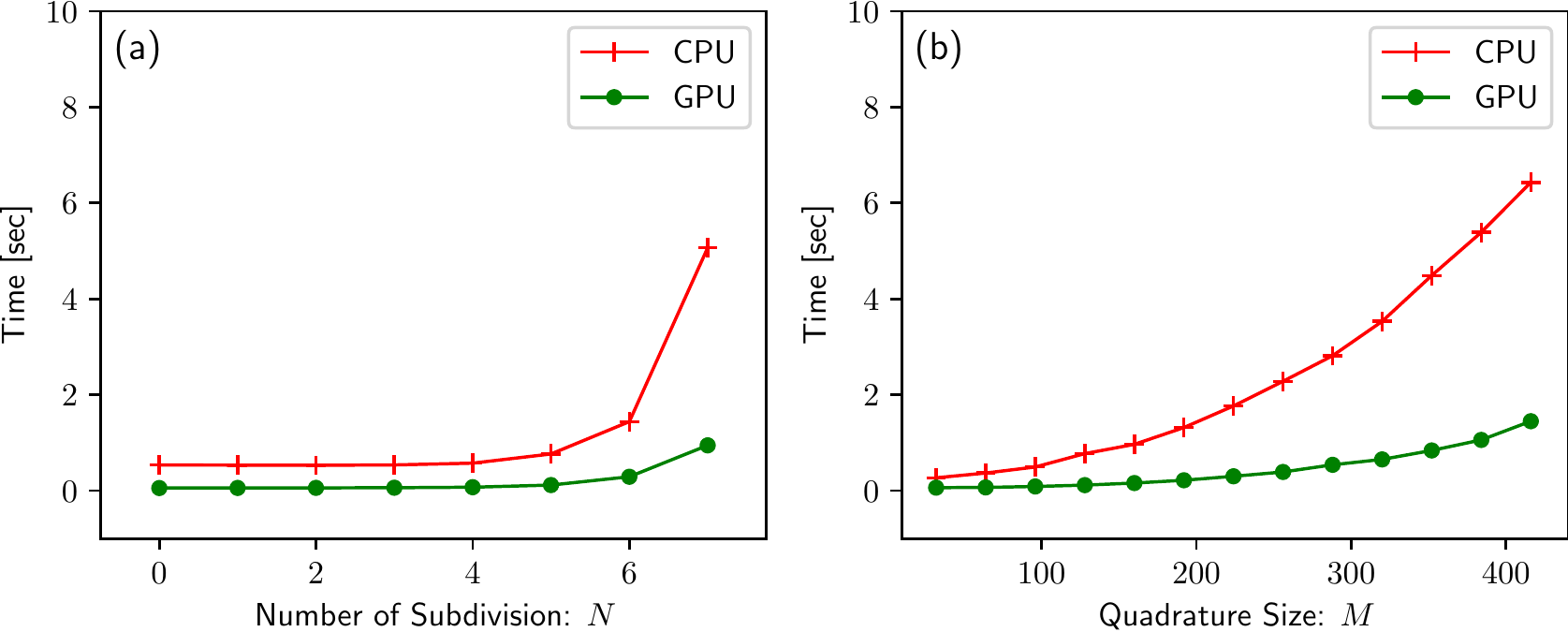}
\end{center}
\caption{\label{fig:time_cpu_gpu}Computation time with parallelization on OpenMP and GPU in single precision : (a) variable $N$ with $(M,\xi_\text{max}) =(128,5)$, and (b) variable $M$ with $(N,\xi_\text{max}) =(5,5)$.} 
\end{figure}
For performance tests, our GPU code is run on NVIDIA RTX A6000 (10\,752 cores and 48 GB global memory) and its computation time is compared against the OpenMP code with the CPU environment described above. Both codes are written in the double precision except for the NUFFT part with the tolerance $\epsilon_\text{rel}=10^{-7}$, and this setting disadvantages consumer-oriented GPUs including RTX A6000, because their peak throughputs in the double precision are 32 times lower than those in the single precision. On the other hand, having repeated calls of the algorithm in mind, the following comparison ignores execution time for memory allocation and warm-up that take longer for the GPU code than for the CPU counterpart. Furthermore, the pointwise operations and reductions on the frequency space are also excluded from both codes, for it can be ideally accelerated by GPUs but are not included in the previous studies \cite{JiWaHuMaNi2019,IcDeTu2023}. Here, Figure\,\ref{fig:time_cpu_gpu}(a) plots computation time versus the number of subdivision $N$ with $M=128$ and $\xi_\text{max} =5$. In both cases, the curves are almost straight for $0\leq N \leq 4$, which implies that the standard FFT dominates the computational costs there, and the growth of time becomes clear for $5\leq N \leq 7$. In particular, our GPU code takes 0.94 sec for $N=7$ while its OpenMP counterpart takes 5.06 sec for the same $N$. As a reference, the direct summation via the exact formula (\ref{eq:simplex32_ft}) takes 10 sec on NVIDIA Tesla V100 for $100\,000$ triangles ($N\approx 6$) and $M=100$ \cite{IcDeTu2023}. On the other hand, Figure\,\ref{fig:time_cpu_gpu}(b) shows computation time versus the quadrature size $M$ with $N=5$ and $\xi_\text{max} =5$, and our GPU code takes 1.44 sec for $M=416$ while the OpenMP version takes 6.42 sec for the same $M$. Again, the GPU code based on the exact formula (\ref{eq:simplex32_ft}) takes roughly 0.1 sec for $80$ triangles ($N=1$) and $M=400$. Since their algorithm has the complexiity $\mathcal{O}(N_S\cdot N_F)$, where $N_S$ is the number of simplices and $N_F$ the number of target points on the frequency space, it is possible to obtain a rough estimate on their computation time for $N=5$ as $ 0.1 \times 20\,480/80 = 25.6$ sec. We thus conclude that our GPU code outperforms the previous study \cite{IcDeTu2023}, and these results are arguably due to the lower complexity $\mathcal{O}(N_F \log N_F + N_P\log^3(1/\epsilon_\text{rel}))$ of the NUFFT-based strategy.
\section{Conclusions}\label{conclusion}
In this paper, we have developed a novel algorithm in the framework of Euclidian harmonic analysis for quantitative comparisons between compact surfaces embedded in $\mathbb{R}^3$. The key idea is to regard those objects as the associated surface measures, which enables to compute their Fourier transforms and then weak distances between them in the sense of the Sobolev norm (\ref{eq:sobolev_2}). This approach allows not only high accuracy due to fast-converging numerical quadrature rules but efficient summations by the Type-1 NUFFT and massive parallelization on GPUs, and the resulted codes outperform the existing methods for similar purposes. As numerical experiments, the $2$-sphere $\mathbb{S}^2$ with the exact Fourier transform  (\ref{eq:ft_nsphere}) is compared with its icosahedral discretization, and it is observed that the simplex-based approximations converge to the smooth object at the quadratic rate up to small truncation. As concluding remarks, we mention three important directions for further developments. \par
First, except for the formulation step in Section \ref{formulate}, the focus of the present work is solely on computational aspects of the suggested algorithm, and neither pointwise convergence on the frequency space nor convergence in the Sobolev norm (\ref{eq:sobolev_2}) is proven. On this point, possible connections between the Hausdorff distance and (\ref{eq:sobolev_2}) can be an interesting topic in numerical analysis of the algorithm, because the Hausdorff distance quantifies differences between compact surfaces in a uniform way with the metric on $\mathbb{R}^3$ and may bound other weaker distances. On the other hand, phenomena called concentration of measures \cite{Ledoux2001} are observed in higher-dimensional spaces and may affect the quantitative comparisons. In the case of the sphere $\mathbb{S}^{d-1}$, its relative surface measure is concentrated around the equator for large $d$, which implies that being close to $\mathbb{S}^{d-1}$ in the sense of a measure-theoretic distance may have far less geometric information. Thus, the suggested algorithm provides some nontirivial questions towards its theoretical studies.  \par
Second, as mentioned in Section \ref{results},  the suggested algorithm involves the user-specified parameter $s$ that controls both the weakness of the Sobolev norm and accuracy of the method, and it should be incorporated into practical applications for validating its effectiveness. For example, it may become alternatives in geometry processing tasks using the Fourier transform in the previous studies \cite{JiWaHuMaNi2019,IcDeTu2023}. In particular, the shape reconstruction in \cite{IcDeTu2023}, which is successful for relatively simple objects with help of the stabilization technique \cite{NiJaJa2021}, is expected to be a good benchmark on to what extent the Sobolev norm (\ref{eq:sobolev_2}) captures geometry of given data. In another direction, our contributions have direct applications to Lagrangian simulations on $\mathbb{S}^2$ \cite{BoKeKrJa2017} where one must carefully check how the current mesh representation deviates from the true smooth surface. Success in those applications will strengthen the claims in the present work and may uncover its limitations at the same time.  \par
Last, the results in the present work motivates us to develop new discretization schemes and their refinement algorithms for more general surfaces. Since the suggested algorithm can reveal quality of geometric approximations in a quantitative manner, it is possible to extend the previous study for periodic planar curves \cite{Koga2022} to higher-dimensional settings. Namely, given an exact parametrization of a smooth surface in $\mathbb{R}^3$, one can approximate the Fourier transform of its surface measure up to the double precision, and then spatial resolution of discretized surfaces are compared to each other by recomputing the Fourier transform from their own data representations.
This procedure is essentially the same as one in Section \ref{results}, and the only difference is that the surface measure of $\mathbb{S}^2$ has the closed form (\ref{eq:ft_nsphere}) and hence there is no need for numerical quadrature rules. For this purpose, it is also important to find spectrally accurate representations for piecewise smooth objects \cite{SaSe2013} and globally smooth but multiply-connected surfaces \cite{ChZhKeLeLuGu2019}. Furthermore, in order to push the order of regularity $s$ towards its critical value, it might be beneficial to consider a periodic analogue of the suggested algorithm for accurate and straightforward computations of Fourier-based distances \cite{RuTu2009}. We hope to address these problems in future works.

\section*{Acknowledgement}
The author would like to thank Shin-ichi Ohta, Keisuke Takasao, and Daisuke Kazukawa for valuable comments on the measure-theoretic aspects of the present work. We are also grateful to the anonymous reviewers for their helpful suggestions, especially those on rigorous aspects of surface measures, their Fourier transforms, and the exact formula (\ref{eq:sphere_sobolev_exact}).

\bibliographystyle{siamplain}
\bibliography{ref}
\end{document}